\documentclass[12pt]{article}
\usepackage[latin1]{inputenc}
\usepackage{amsmath,amsthm,amssymb}

\newtheorem{thm}{Theorem}[section]
 
 \newtheorem{lem}[thm]{Lemma}

\def\endeq{\end{equation}}

\def\eq#1{(\ref{#1})}

\newcommand{\neweq}[1]{\begin{equation}\label{#1}}

\def\ep{\varepsilon}

\def\phi{\varphi}
\def\RR{\mathbb R}

\def\di{\displaystyle}
\def\ri{\rightarrow}

\def\incep{\left\{\begin{array}{cl} }
 \def\termin{\end{array}\right. }
\def\2af{2^*_\alpha}

\def\proof{{\it Proof.}\ }
\def\di{\displaystyle}
\def\ri{\rightarrow}

\title{\sc Asymptotics for singular solutions of quasilinear elliptic equations with absorption term}
\author{Du\v{s}an Repov\v{s}$\,^{a,b}$\\
 \small $^a\,$Faculty of Mathematics and
Physics, University of Ljubljana,\\ \small Jadranska  19,  P. O. Box 2964, 1001 Ljubljana, Slovenia\\
\small $^b\,$Faculty of Education, University of Ljubljana,\\ \small Kardeljeva plo\v{s}\v{c}ad 16, 1000 Ljubljana, Slovenia\\
 \small E-mail: {\tt dusan.repovs@guest.arnes.si}\\}

\date{}

\begin{document}

\maketitle

\begin{abstract}
We are concerned with the asymptotic analysis of positive blow-up boundary  solutions for a class of quasilinear elliptic equations with absorption term. By means of the Karamata theory we establish the first two terms in the expansion of the singular solution near the boundary. Our analysis includes large classes of nonlinearities of Keller-Osserman type.
\end{abstract}

\section{Introduction and the main result}
Let $\Omega\subset\RR^N$ ($N\geq 2$) be a bounded domain with $C^2$ boundary. Throughout this paper we assume that $1<p<\infty$, $a:\overline\Omega\rightarrow (0,\infty)$ is a H\"older potential,
and $f:[0,\infty)\ri [0,\infty)$ is a $C^1$ function.

We are concerned with the study of solutions $u\in W^{1,p}_{\rm loc}(\Omega)\cap C^{1,\mu}(\Omega)$ of the following quasilinear elliptic problem
\begin{equation}\label{P}
\left\{\begin{array}{lll}
&\di \Delta_pu=a(x)f(u)&\qquad\mbox{in}\ \Omega\\
&\di u(x)\ri +\infty&\qquad\mbox{as }\mbox{dist}(x,\partial\Omega)\ri 0\\
&\di u>0&\qquad\mbox{in}\ \Omega\,.\end{array}\right.
\end{equation}
Under appropriate assumptions, the existence of a solution for problem \eq{P} has been proved in \cite{moha1}.
Our objective in this paper is to establish the first two terms of the boundary blow-up rate for solutions of \eq{P}, under appropriate conditions on the nonlinearity $f$ and the variable potential $a$.

This problem can be regarded as a model of a steady-state single species inhabiting 
$\Omega$, so $u(x)$ stands for the population density. In fact, if $f(u)=u^q$ ($q>p-1$), problem \eq{P} is a basic population model and it is also related to
some prescribed curvature problems in Riemannian geometry.  We refer the reader to Li,  Pang and Wang \cite{lipang} for a study of problem \eq{P}
in the case of multiply connected domains and subject to mixed boundary conditions.

The study of singular problems with blow-up on the boundary was initiated in the case $p=2$, $a\equiv 1$, and $f(u)=\exp(u)$ by
Bieberbach \cite{bi} (if $N=2$) and Rademacher (if ($N=3$). Problems of this type
arise in Riemannian geometry, namely if a Riemannian metric of the form
$|ds|^2=\exp(2u(x))|dx|^2$ has constant Gaussian curvature $-c^2$
then $\Delta u=c^2 \exp(2u)$. Such problems also appear in the theory of  automorphic functions, Riemann surfaces, as well as  in the theory of the electric potential in a glowing hollow metal body. Lazer
and McKenna \cite{lm} extended the results of Bieberbach and
Rademacher for bounded domains in ${\mathbb R}^N$ satisfying a
uniform external sphere condition and for exponential-type nonlinearities.
An important development is due to Keller \cite{ke} and
Osserman \cite{os}, who established a necessary and sufficient condition for problem \eq{P} to have
a solution, provided that $p=2$, $a\equiv 1$, and $f$ is an increasing nonlinearity.
In a celebrated paper connected with the Yamabe problem, Loewner and Nirenberg \cite{ln}
linked the uniqueness of the blow-up solution to
the growth rate at the boundary.
Motivated by certain geometric problems, they
established the uniqueness for the case $f(u)=u^{(N+2)/(N-2)}$, $N>2$.
For related results we refer the reader to Bandle and Marcus \cite{bm}, Bandle, Moroz and Reichel \cite{bamr}, L\'opez-G\'omez \cite{lopez}, Marcus and V\'eron \cite{mv1}, Mohammed \cite{moha},
Repov\v{s} \cite{repovs}, etc. The case of nonmonotone nonlinearities was studied by Dumont, Dupaigne, Goubet and R\u adulescu \cite{ddgr}.

In order to describe our main result, we need to recall  some basic notions and properties in the theory of
  functions with regular variation at infinity and of functions belonging to the Karamata class. We point out that Karamata \cite{karamata}  introduced this theory in relation
  to
  Tauberian theorems. This theory was then applied to the analytic number theory, analytic functions, Abelian theorems, and probability theory (see Feller \cite{feller}). We refer
  the reader to the works by Bingham, Goldie, and  Teugels \cite{bgt} and Seneta \cite{seneta} for details and related results. The combined use of the regular variation theory and the Karamata theory has been introduced by C\^{\i}rstea and R\u
adulescu \cite{crasun,crasas,crasan,crtams} in the study of various qualitative and asymptotic properties of solutions of nonlinear partial differential equations. In particular, this setting
becomes a powerful tool in describing the asymptotic behavior of solutions for large classes of nonlinear elliptic equations, including singular solutions with blow-up boundary and stationary problems with either degenerate or singular nonlinearity.

\smallskip We say that a positive
measurable function $f$ defined on some interval $[B,\infty)$ is  regularly varying at infinity
with index  $q\in \RR$  if
for all $\xi>0$
$$ \lim_{u\to \infty}f(\xi u)/f(u)=\xi^q. $$
When the index of regular variation $q$ is zero, we say that the
function is slowly varying.

If $RV_q$ denotes the class of functions with regular variation with index $q$
then the function
 $f(u)=u^q$ belongs to $RV_q$. The functions $\ln (1+u)$,
$\ln \ln(e+u)$, ${\rm exp}\,\{(\ln u)^\alpha\}$, $\alpha\in (0,1)$
vary slowly, as well as any measurable function with positive limit at infinity.
Using the definition of $RV_q$, a straightforward computation shows that if $p>1$ and $f\in RV_q$ with $q>p$ is continuous and increasing on $[B,\infty)$ then its anti-derivative $F(t):=\int_B^t f(s)ds$ satisfies $F\in RV_{q+1}$, and hence $F^{-1/p}\in RV_{-(q+1)/p}$. According to \cite{crasun} (see also \cite{moha}), we deduce that $F^{-1/p}\in L^1(B,\infty)$, that is, $f$ satisfies the Keller-Osserman condition
\begin{equation}\label{kelose}
\int^\infty \left[ F(t)\right]^{-1/p}<\infty\,.\end{equation}

An important subclass of $RV_q$ contains the functions $f$ such that $u^{-q}f(u)$ is a renormalized slowly varying function. More precisely, we denote by $NRV_q$  the set of functions $f$ having the form $f(u)=Au^q\exp (\int_B^u\phi (t)/tdt)$ for all $u\geq B>0$, where $A$ is a positive constant and $\phi\in C[B,\infty)$ satisfies $\lim_{t\ri\infty}\phi (t)=0$. Then, by the Karamata representation theorem (see \cite{bgt}), we have $NRV_q\subset RV_q$.

Next, we denote by
${\cal K}$ the class of all positive, increasing
$C^1$-functions $k$
defined on $(0,\nu)$, for some $\nu>0$, which satisfy
$\lim_{t\to 0^+}\di\left(\frac{K(t)}{k(t)}\right)^{(i)}:=\ell_i$ for
$i\in\{0,1\}$, where $K(t)= \int_0^t
k(s)\,ds$. A straightforward computation shows that
$\ell_0=0$ and $\ell_1\in [0,1]$, for all $k\in {\cal K}$.

Let ${\cal K}_{0,1}$ denote the set of all functions
$k\in{\cal K}$ satisfying
$$\lim_{t\searrow 0}t^{-1}
\left[\left(K(t)/k(t)\right)'-\ell_1\right]:=L_1\in \RR\,.$$.

We study problem \eq{P} provided that the nonlinear term $f$ satisfies
\begin{equation}\label{f1}
f\in C^1[0,\infty),\quad f(0)=0,\quad f>0 \mbox{ and } f \mbox{ is increasing on $(0,\infty)$}.\end{equation}

We now describe the growth of $f$ at infinity. We assume that $f\in NRV_{\sigma +1}$ for some $\sigma >p-2$. This means that $f$ can be written as
$$f(u)=A_0u^{\sigma +1}\exp\left(\int_B^u \phi (t)/t\,dt\right),$$
for some $A_0>0$, where $\phi\in C^1[B,\infty)$ and $\lim_{t\ri\infty}\phi (t)=0$. Moreover, we assume that there is some $\frac{\sigma+2}{p}-1<\alpha <\sigma +2$ such that
\begin{equation}\label{f2}
\lim_{t\ri\infty}\frac{t\phi '(t)}{\phi (t)}=-\alpha\,.\end{equation}

We also assume that $a:\overline\Omega\ri (0,\infty)$ satisfies $a\in C^{0,\mu}(\overline\Omega)$ for some $0<\mu<1$ and $k\in {\cal K}_{0,1}$,
\begin{equation}\label{aa}
a(x)=k^p(d(x))\left( 1+Ad(x)+o(d(x))\right)\qquad\mbox{as}\ d(x)\ri 0,
\end{equation}
where $A>0$ and $d(x):=\mbox{dist}\, (x,\partial\Omega)$.

For any $x\in\Omega$ near the boundary of $\Omega$ we denote by $\overline x\in\partial\Omega$ the unique point such that $d(x)=|x-\overline x|$. We also denote by ${\cal H}(\overline x)$ the mean curvature of $\partial\Omega$ at the point $\overline x$.

Our main result extends to
a quasilinear setting the results
given in \cite{crasas}, \cite{moha}, and \cite{zhang}. Our asymptotic development also relies on the geometry of the domain, as developed by Bandle and Marcus \cite{bmdie}.

\begin{thm}\label{teo}
Assume that 
$f\in NRV_{\sigma +1}$ ($\sigma >p-2$) satisfies hypotheses \eq{f1} and \eq{f2}.
Suppose 
that $a\in C^{0,\mu}(\overline\Omega)$ satisfies condition \eq{aa}. Then any solution of problem \eq{P} satisfies
$$u(x)=\xi_0 h(K(d(x)))\left( 1+ C_1d(x)+C_2{\cal H}(\overline x)d(x)+o(d(x))\right)\quad\mbox{as $d(x)\ri 0$},$$
where $h$ is uniquely defined by
$$\left(\frac{p-1}{p}\right)^{1/p}\int_{h(t)}^\infty \left( F(t)\right)^{-1/p}dt=t$$
and
$$\xi_0=\left[(p-1)\,\frac{p+\ell_1(\sigma+2-p)}{\sigma+2}\right]^{1/(\sigma+2-p)},$$
$$C_1=\frac{L_1(\sigma+2-p)-A(p+(\sigma+2-p)\ell_1)}{\sigma[\ell_1(\sigma+2-p)+p]},$$
$$C_2=\frac{\ell_1(N-1)(\sigma+2-p)}{\ell_1(\sigma+2-p)+(\sigma+1)(\sigma+2)-p}\,.$$
\end{thm}

\section{Auxiliary results}
The proof of the main result strongly relies on the maximum principle for quasilinear equations in the following form. We refer the reader to \cite{tolk} for a detailed proof and related results.

\begin{lem}\label{l1}
Let $\Omega$ be a bounded domain in $\RR^N$ with smooth boundary. Assume
that $V_1$ and $V_2$ are continuous functions on $\Omega$ such that $V_1\in L^\infty (\Omega)$ and $V_2>0$. Let $u_1,\, u_2\in W^{1,p}(\Omega)$ be positive  functions such that
\begin{equation}\label{l11}\begin{array}{ll}
&\di \Delta_p u_1+V_1u_1^{p-1}+V_2f(u_1)\leq 0\leq \\ &\di\Delta_p u_2+V_1u_2^{p-1}+V_2f(u_2)\qquad\mbox{in ${\mathcal D}'(\Omega)$}\end{array}
\end{equation}
and
\begin{equation}\label{l12}
\limsup_{x\ri\partial \Omega}(u_2(x)-u_1(x))\leq 0,
\end{equation}
where $f$ is continuous on $[0,\infty)$ such that the mapping $f(t)/t^{p-1}$ is increasing for $\inf_\Omega (u_1,u_2) <t<
\sup_\Omega (u_1,u_2)$.

Then $u_1\geq u_2$ in $\Omega$.
\end{lem}

The proof of Lemma \ref{l1}  relies on some ideas introduced by Benguria, Brezis and Lieb \cite{bbl} (see also Marcus and V\'eron
\cite[Lemma 1.1]{mv1}, C\^{\i}rstea and R\u adulescu \cite[Lemma 1]{crhouston}, and Du and Guo \cite{duguo}).

Our growth rate of $f$ expressed by
the assumptions $f\in NRV_{\sigma +1}$ and $\sigma >p-2$ implies
that $f$ satisfies the Keller-Osserman condition \eq{kelose} and
$$\lim_{t\ri\infty}\frac{tf(t)}{F(t)}=\sigma +2\,.$$
Next, we set
$${\cal F}(t):=\left(\frac{p-1}{p}\right)^{1/p}\int_{t}^\infty ( F(x))^{-1/p}dx\,.$$
Since
$${\cal F}'(t)=-\left(\frac{p-1}{p}\right)^{1/p}( F(t))^{-1/p}\,,$$
we deduce that
$$\lim_{t\ri\infty}\frac{t{\cal F}'(t)}{{\cal F}(t)}=-\frac{\sigma +2}{p}-1$$
and
$$\lim_{t\ri\infty}\frac{F(t)^{(p-1)/p}}{f(t){\cal F}(t)}=\frac 1p\, \left(\frac{p}{p-1}\right)^{1/p}\left(1-\frac{p}{\sigma+2}\right)\,.$$

These estimates enable us to deduce the following auxiliary result.

\begin{lem}\label{l2}
Under the assumptions of Theorem \ref{teo}, the following properties hold true:

(i) $\lim_{t\ri\infty}\frac{\frac{tf'(t)}{f(t)}-\sigma-1}{{\cal F}(t)}=\lim_{t\ri\infty}\frac{\frac{F(t)}{tf(t)}-\frac{1}{\sigma +2}}{{\cal F}(t)}=0$;

(ii) $\lim_{t\ri\infty}\frac{\left(\frac{p}{p-1}\right)^{(p-1)/p}\frac{(F(t))^{(p-1)/p}}{f(t){\cal F}(t)}\,-\,\frac{\sigma+2-p}{(p-1)(\sigma+2)}}{{\cal F}(t)}=0$;

(iii) $\lim_{t\ri\infty}\frac{\frac{f(at)}{a^{p-1}f(t)}-a^{\sigma+2-p}}{{\cal F}(t)}=0$, for all $a>0$.
\end{lem}

\proof The proofs of (i) and (ii) follow directly by the previous considerations about $f$, $F$, and ${\cal F}$.

(iii) If $a=1$ the property is obvious. Let us now assume that $a\not=1$. We have
$$\frac{f(at)}{a^{p-1}f(t)}-a^{\sigma+2-p}=a^{\sigma+2-p}\left[\exp\left(\int_t^{at}\frac{\phi(x)}{x}dx
\right)-1\right]\,.$$
Our hypotheses on $\phi$ imply that
$$\lim_{t\ri\infty}\frac{\phi(tx)}{x}=0\qquad\mbox{and}\qquad
\lim_{t\ri\infty}\frac{\phi(tx)}{x\phi (x)}=x^{-\alpha-1}\,,$$
uniformly for either $x\in [a,1]$ or $x\in [1,a]$. This implies that
$$\lim_{t\ri\infty}\int_t^{at}\frac{\phi(x)}{x}dx=\int_1^{a}\frac{\phi(tx)}{x}dx=0$$
and
$$\lim_{t\ri\infty}\int_1^{a}\frac{\phi(tx)}{x\phi(t)}dx=
\lim_{t\ri\infty}\int_1^{a}x^{-\alpha-1}dx=\alpha{-1}(1-a^{-\alpha}).$$
We conclude that
$$\begin{array}{ll} \di \frac{f(at)}{a^{p-1}f(t)}-a^{\sigma+2-p} &\di = a^{\sigma+2-p}\lim_{t\ri\infty}
\frac{\int_1^a\frac{\phi (tx)}{x}dx}{{\cal F}(t)}\\
&\di =a^{\sigma+2-p}\lim_{t\ri\infty}\frac{\phi(t)}{{\cal F}(t)}\,
\lim_{t\ri\infty}\int_1^a
\frac{\phi (tx)}{x\phi (t)}\,dx=0.\end{array}$$
This completes the proof. \qed

We conclude this section with some properties of the function $h$ that describes the blow-up rate of solutions of problem \eq{P} in the statement of Theorem \ref{teo}.

\begin{lem}\label{l3}
Assume that 
the hypotheses of Theorem \ref{teo} are fulfilled and let $h:(0,\infty)\ri (0,\infty)$
be the function defined implicitly by
$$\left(\frac{p-1}{p}\right)^{1/p}\int_{h(t)}^\infty \left( F(t)\right)^{-1/p}dt=t.$$
Then the following properties hold:

(i) $\lim_{t\searrow 0}th'(t)/h(t)=-p/(\sigma +2-p)$;

(ii) $\lim_{t\searrow 0}h'(t)/(th''(t))=-(\sigma +2-p)/(\sigma+2)$;

(iii) $\lim_{t\searrow 0}h(t)/(t^2h''(t))=(\sigma +2-p)^2/[p(\sigma+2)]$;

(iv) $\lim_{t\searrow 0}\left(\frac{h'(t)}{th''(t)}+\frac{\sigma+2-p}{\sigma+2}\right)/t=0$;

(v) for all $k\in {\cal K}_{0,1}$,
$$\lim_{t\searrow 0}t^{-1}\left(1+\frac{k'(t)K(t)}{k^2(t)}\cdot\frac{h'(K(t))}{K(t)h''(K(t))}-
\frac{1}{p-1}\cdot\frac{f(\xi_0h((K(t))))}{\xi_0^{p-1}f(h(K(t)))}\right)=
\frac{(\sigma+2-p)L_1}{\sigma+2}\,.$$
\end{lem}

\proof We first observe that $\lim_{t\ri 0}h(t)=+\infty$ and $h'(t)=-p^{1/p}(p-1)^{-1/p}F(t)^{1/p}$.

(i) We have
$$\lim_{t\searrow 0}\frac{th'(t)}{h(t)}=-\lim_{s\ri +\infty}\frac{(F(s))^{1/p}\int_s^\infty (F(v))^{-1/p}dv}{s}=-\frac{p}{\sigma+2-p}\,.$$

(ii) A straightforward computation shows that for all $t>0$, $$h''(t)=(p-1)^{-2/p}p^{(2-p)/p}f(h(t))(F(h(t)))^{(2-p)/p}\,.$$
Therefore
$$\lim_{t\searrow 0}\frac{h'(t)}{th''(t)}=-p^{1/p}(p-1)^{(p-1)/p}\cdot\lim_{s\ri +\infty}\frac{(F(s))^{1/p}}{f(s){\cal F}(s)}=-\frac{\sigma +2-p}{\sigma +2}\,. $$

(iii) We have
$$\lim_{t\searrow 0}\frac{h(t)}{t^2h''(t)}=\lim_{t\searrow 0}\frac{h(t)}{th'(t)}\cdot
\lim_{t\searrow 0}\frac{h'(t)}{th''(t)}=\frac{(\sigma +2-p)^2}{p(\sigma+2)}\,.$$

(iv) The proof follows by combining the previous results.

(v) The proof follows after combining  Lemma \ref{l2} with the previous results.
\qed

\section{Proof of Theorem \ref{teo}}
For fixed $\eta>0$ small enough, we define
$$\Omega_\eta:=\{x\in\Omega:\ 0<d(x)<\eta\}\,.$$

For any $x\in\Omega$, we set $r=d(x)=|x-\overline x|$. Define
$$S_1(r)=r^{-1}\left(1+\frac{k'(r)K(r)}{k^2(r)}\cdot\frac{h'(d(r))}{K(r)h''(K(r))}-
\frac{1}{p-1}\cdot\frac{f(\xi_0h(K(r)))}{\xi_0^{p-1}f(h(K(r)))}\right)\,.$$
Then, by Lemma \ref{l3}, we have $\lim_{r\searrow 0}S_1(r)=L_1(\sigma+2-p)/(\sigma+2)$.

Fix $\ep>0$ small enough. Since $\Omega$ has smooth boundary,  there exists $\delta =\delta(\Omega)>0$ such that $d\in C^2(\overline\Omega_\delta)$ and for all $x\in\Omega_\delta$,
$|\nabla d(x)|=1$. Set, for all $x\in\Omega_\delta$,
$$z_\pm (x)=\xi_0 h(K(d(x))\left(1+(C_1\pm\ep)d(x)+C_2{\cal H}(\overline x)d(x)\right)\,.$$
Then, by the mean value theorem, there exists $\lambda_\pm \in (0,1)$ depending on $x$ such that for all $x\in\Omega_\delta$,
$$f(z(x))=f(\xi_0(h(K(d(x)))))+\xi_0h(K(d(x)))f'(h_\pm(d(x)))
\left((C_1\pm\ep)d(x)+C_2{\cal H}(\overline x)d(x) \right),$$
where
$$h_\pm (d(x))=\xi_0(h(K(d(x))))\left(1+\lambda_\pm((C_1\pm\ep)d(x)+C_2{\cal H}(\overline x)d(x))\right)\,.$$

Define the mapping
$$\begin{array}{ll}\di
S_{2\pm}(r)&\di =(C_1\pm\ep)\left[1+\frac{h'(K(r))}{K(r)h''(K(r))}\left(\frac{K(r)k'(r)}{k^2(r)}+
\frac{2K(r)}{rk(r)} \right)\right]\\
&\di - \frac{C_1\pm\ep}{p-1}\, \frac{f'(h_\pm(K(r)))}{f'(h(K(r)))}\,
\frac{h(K(r))f'(h(K(r)))}{\xi_0^{p-2}f(h(K(r)))}\\
&\di -\frac{1}{p-1}\, (A\mp\ep )\frac{f(\xi_0h(K(r)))}{\xi_0^{p-1}f(h(K(r)))},
\end{array}$$
where $0<\eta<\min{1,p-2}$.
Using Lemma \ref{l3} we deduce that the asymptotic behavior of $S_{2\pm}$ near the origin is given by
$$\begin{array}{ll} \di \lim_{r\ri 0}S_{2\pm}(r)&\di =-
\left( C_1\frac{\ell_1(\sigma+2-p)(\sigma+2)+p}{\sigma+2}+A\,
\frac{p+\ell_1(\sigma+2-p)}{\sigma+2}\right)\\
&\di\mp \ep\left(\frac{\ell_1(\sigma+2-p)(\sigma+2)+p}{\sigma+2}+
\eta\,\frac{p+\ell_1(\sigma+2-p)}{\sigma+2}\right).\end{array}$$

We also define the mappings
$$\begin{array}{ll}\di
S_{3}(x)&\di =C_2{\cal H}(\overline x)\left[1+\frac{h'(K(r))}{K(r)h''(K(r))}\left(\frac{K(r)k'(r)}{k^2(r)}+
\frac{2K(r)}{rk(r)}\right) \right]\\
&\di - \frac{{\cal H}(\overline x)}{p-1}\, \frac{f'(h_\pm(K(r)))}{f'(h(K(r)))}\,
\frac{h(K(r))f'(h(K(r)))}{\xi_0^{p-2}f(h(K(r)))}\\
&\di -(N-1){\cal H}(\overline x)\,\frac{h'(K(r))}{K(r)h''(K(r))}\,\frac{K(r)}{rk(r)}\,,
\end{array}$$
$$\begin{array}{ll}\di
S_{4\pm}(x)&\di =r\,\frac{h'(K(r))}{K(r)h''(K(r))}\, (C_1\pm\ep+C_2{\cal H}(\overline x))\Delta d(x)\\
&\di +(C_1\pm\ep+C_2{\cal H}(\overline x))\frac{h(K(r))}{K^2(r)h''(K(r))}\,\frac{K^2(r)}{rk^2(r)}\Delta d(x)\\
&\di - (A\mp\eta\ep)(C_1\pm\ep+C_2{\cal H}(\overline x))r\,\frac{f'(h_\pm (K(r)))}{f'(h(K(r)))}\,\frac{h(K(r))f'(h(K(r)))}{\xi_0^{p-2}f(h(K(r)))}\,.
\end{array}$$
Applying again Lemma \ref{l3} we deduce that
$$\lim_{d(x)\ri 0}S_{3}(x)=\lim_{d(x)\ri 0}S_{4\pm}(x)=0.$$
Therefore
$$\begin{array}{ll}
&\di\lim_{d(x)\ri 0}(S_{1}(r)+S_{2\pm}(r)+S_{3}(x)+S_{4\pm}(x)) =\\ &\di
\mp\frac{\ep}{\sigma+2}\left[p+\ell_1(\sigma+2-p)(\sigma+2)+
\eta(p+\ell_1(\sigma+2-p))\right]\,.\end{array}$$
Finally, we define
$$
S_{5\pm}(x)\di =\left|\left[\left(1+(C_1\pm\ep )r+C_2{\cal H}(\overline x)r\right)+\left((C_1\pm\ep)+C_2{\cal H}(\overline x))\right)\frac{K(r)}{k(r)}\,\frac{h(K(r))}{K(r)h'(K(r))}\right]\nabla d(x)
 \right|
\,.$$
We observe that our hypotheses imply
$$\lim_{d(x)\ri 0}S_{5\pm}(x)=0\,.$$

Our hypotheses imply that there are positive numbers $\delta_{1\ep}$ and $\delta_{2\ep}$ such that $0\leq K(t)\leq 2\delta_{1\ep}$ for all $t\in (0,2\delta_{2\ep})$ and
for all $x\in\Omega_{2\delta_{1\ep}}$,
$$k^p(d(x))(1+(A-\eta\ep)d(x))\leq a(x)\leq
k^p(d(x))(1+(A+\eta\ep)d(x))\,.$$
At the same time, restricting eventually $\delta_{1\ep}$ and $\delta_{2\ep}$, we can assume that for all $x\in \Omega_{2\delta_{1\ep}}$ with $|x-\overline x|<2\delta_{2\ep}$,
$$ S_{1}(r)+S_{2+}(r)+S_{3}(x)+S_{4+}(x)\leq 0\leq  S_{1}(r)+S_{2-}(r)+S_{3}(x)+S_{4-}(x) .$$

Next, for some fixed $\rho\in (0,2\delta_{1\ep})$, we define
$d_1(x)=d(x)-\rho$, $d_2(x)=d(x)+\rho$, and
$$\Omega_\rho^-=\{x\in\Omega;\ \rho <d(x)<2\delta_{1\ep}\}\qquad
\Omega_\rho^+=\{x\in\Omega;\ d(x)<2\delta_{1\ep}-\rho\}\,.$$

Set
$$\overline u_\ep(x)=\xi_0 h(K(d_1(x)))\left( 1+ (C_1+\ep)d_1(x)+C_2{\cal H}(\overline x)d_1(x)\right)\qquad x\in\Omega_\rho^-$$
and
$$\underline u_\ep(x)=\xi_0 h(K(d_2(x)))\left( 1+ (C_1-\ep)d_2(x)+C_2{\cal H}(\overline x)d_2(x)\right)\qquad x\in\Omega_\rho^+\,.$$

Our main purpose in what follows is to show that $\overline u_\ep$ is a supersolution of the equation \eq{P} in $\Omega_\rho^-$
and $\underline u_\ep$ is a subsolution of \eq{P} in $\Omega_\rho^+$. We first observe that the mean value theorem implies
$$\begin{array}{ll} \di f(\overline u_\ep) &\di = f(\xi_0h(K(d_1(x))))\\
&\di + \xi_0h(K(d_1(x)))f'(h_+(d_1(x)))[(C_1+\ep)d_1(x)+C_2{\cal H}(\overline x)d_1(x)]\end{array}$$
for all $x\in\Omega_\rho^-$,
where, for some $\zeta\in (0,1)$ depending on $x$,
$$h_+(d_1(x))=\xi_0 h(K(d_1(x)))[1+\zeta (C_1+\ep)d_1(x)+C_2{\cal H}(\overline x)d_1(x)]\,.$$
Combining these results, we deduce that for all $x\in\Omega_\rho^-$,
$$\begin{array}{ll}
&\di \Delta_p\overline u_\ep (x)-k^p(d_1(x))\, (1+(A-\ep )d_1(x)f(\overline u_\ep))=\\
&\di  (p-1)\xi_0^{p-1}\, k^p(d_1(x))d_1(x)\,|h'(K(d_1(x)))|^{p-2}\, h''(K(d_1(x))) \cdot\\
&\di S_{5+}(x)\left(S_1(r)+S_{2+}(r)+S_{3}(x)+S_{4+}(x)\right)\leq 0,
\end{array}$$
where $r=d_1(x)+\rho$.

We now deduce uniform estimates for the solution of problem \eq{P} in terms of $\overline u_\ep$ and $\underline u_\ep$. For this purpose we follow the method introduced in \cite{crasun}. Assume
that $u$ is an arbitrary solution of problem \eq{P}. Thus, for all $x\in \partial\Omega_\rho^-$,
$$u(x)\leq \overline u_\ep (x)+M_1(\delta_{1\ep}),\qquad\mbox{where
$M_1(\delta_{1\ep})=\max_{d(x)\geq\delta_{1\ep}}u(x)$}.$$
Thus, by the maximum principle,
\begin{equation}\label{equ11}u(x)\leq \overline u_\ep (x)+M_1(\delta_{1\ep}),\qquad\mbox{for all $x\in \Omega_\rho^-$}.\end{equation}

Next, since the function $h$ is decreasing, we have for all $x\in\Omega$ with $d(x)=
2\delta_{1\ep}-\rho$,
$$\underline u_\ep (x)\leq\xi_0h(K(2\delta_{1\ep})):=M_2(\delta_{1\ep}).$$
The maximum principle implies that
\begin{equation}\label{equ12}\underline u_\ep (x)\leq u(x)+M_2(\delta_{1\ep})\qquad\mbox{for all $x\in\Omega_\rho^+$.}\end{equation}
 Taking $\rho\ri 0$ in relations \eq{equ11} and \eq{equ12} we obtain, for all $x\in \Omega_\rho^-\cap \Omega_\rho^+$,
 $$\begin{array}{ll}
 &\di 1+(C_1-\ep)d(x)+C_2{\cal H}(\overline x)d(x)-\frac{M_2(\delta_{1\ep})}{\xi_0h(K(d(x)))}\leq\frac{u(x)}{\xi_0h(K(d(x)))}\leq\\
 &\di 1+(C_1+\ep)d(x)+C_2{\cal H}(\overline x)d(x)+\frac{M_2(\delta_{1\ep})}{\xi_0h(K(d(x)))}\,.\end{array}$$
 This implies that
 $$\begin{array}{ll}
  C_1-\ep +C_2{\cal H}(\overline x)&\di\leq\liminf_{d(x)\ri 0}\frac{1}{d(x)}\left(
 \frac{u(x)}{\xi_0h(K(d(x)))}-1\right)\\
 &\di\leq \limsup_{d(x)\ri 0}\frac{1}{d(x)}\left(
 \frac{u(x)}{\xi_0h(K(d(x)))}-1\right)\\ &\di\leq C_1+\ep +C_2{\cal H}(\overline x).
 \end{array}$$
 Taking now $\ep\ri 0$ we conclude that
 $$u(x)=\xi_0 h(K(d(x)))\left( 1+ C_1d(x)+C_2{\cal H}(\overline x)d(x)+o(d(x))\right)\quad\mbox{as $d(x)\ri 0$}.$$
 This completes the proof.\qed

\medskip
{\bf Acknowledgments}. The author acknowledges the support by the Slovenian Research Agency
grants P1-0292-0101, J1-2057-0101 and J1-4144-0101.

\end{document}